\documentclass{article}
\usepackage{amssymb,amsthm,amsmath,amsfonts}

\theoremstyle{plain}
\begingroup

\endgroup

\def\be{\begin{equation}}
\def\ee{\end{equation}}
\def\ba{\begin{eqnarray}}
\def\ea{\end{eqnarray}}

\newcommand{\R}{\mathbb R}

\newcommand{\man}{\mathcal{M}}
\newcommand{\manb}{\mathcal{N}}
\newcommand{\tanv}{\mathcal{V}}
\newcommand{\metric}{g}
\newcommand{\energy}{\mathcal{E}}
\newcommand{\vect}{v}
\newcommand{\tr}{\mathrm{Tr}}
\newcommand{\Per}{\mathrm{Per}}
\newcommand{\id}{\mathrm{Id}}

\newcommand{\Fs}{\mathbf{F}_s}
\newcommand{\Fl}{\mathbf{F}_l}
\newcommand{\gradperp}{\nabla_\perp}
\newcommand{\lap}{\triangle}
\newcommand{\vel}{\mathbf{u}}
\newcommand{\Vel}{\mathbf{U}}
\newcommand{\vell}{\mathbf{v}}
\newcommand{\velll}{\mathbf{w}}
\newcommand{\Vell}{\mathbf{V}}
\newcommand{\Velll}{\mathbf{W}}
\newcommand{\velf}{\mathbf{f}}
\newcommand{\vpa}{\varphi}
\newcommand{\vpb}{\Phi}

\newcommand{\ei}{\mathbf{i}}
\newcommand{\ej}{\mathbf{j}}
\newcommand{\ek}{\mathbf{k}}
\newcommand{\normal}{\mathbf{n}}
\newcommand{\qref}[1]{(\ref{#1})}

\theoremstyle{definition}

\begin{document}

\title{Gradient flow structure for domain relaxation in Langmuir films}

\author{ Mahir Had\v{z}i\'{c}\textsuperscript{1} and Govind Menon\textsuperscript{2}}
\date{}
\maketitle

\medskip
\noindent
\footnotetext[1]
{Institut f\"ur Mathematik, Universit\"at Z\"urich, Winterthurerstrasse 190, CH-8057 Z\"urich, Switzerland.
Email: mahir.hadzic@math.uzh.ch}
\footnotetext[2]{ Division of Applied Mathematics, Box F, Brown University, Providence, RI 02912.
Email: menon@dam.brown.edu}

\begin{abstract}
We describe a gradient flow structure for the inviscid Langmuir layer Stokesian subfluid model introduced recently by Alexander {\em et al}~\cite{AlBe}.  
\end{abstract}

\section{Introduction}
In a recent article Alexander {\em et al}~\cite{AlBe} introduced a model for the evolution of molecularly thin Langmuir layers on the surface of a subfluid. The Langmuir layer consists of domains of different phases, and the simplest model is a study of  layers with only two phases. For a quiescent subfluid, the evolution of the Langmuir layer is driven by a line tension at the interface between phases  and damped  by a Stokes flow in the subfluid. The resulting free boundary problem is called the inviscid Langmuir layer Stokesian subfluid (ILLSS) model.  

In this brief note we comment on the mathematical structure of the ILLSS model. Alexander {\em et al} showed that the perimeter of the interface determines a natural energy of the system which is dissipated by the Stokesian subfluid. We show that this energy-energy dissipation relation arises from a {\em gradient flow\/}. 

This observation is based on a comparison between the ILLSS model and the Hele-Shaw flow. Alexander {\em et al} showed that the ILLSS model admits a boundary integral formulation similar to that of the Hele-Shaw flow, and  it is known that the Hele-Shaw flow has a gradient structure.  The first such formulation is attributed to Fife in~\cite{Almgren}. We follow a distinct approach due to Otto~\cite{Otto}. The ingredients of a gradient flow are a manifold (the phase space) with a metric and an energy functional.
We find that as in the Hele-Shaw flow, the energy in the ILLSS model is the perimeter of domains. However, it differs from the Hele-Shaw flow in the choice of manifold and metric. The main heuristic idea is that the metric must reflect the dissipation, and here it occurs in the Stokes flow in the subfluid. 

There is also a striking difference with the Hele-Shaw flow. An interesting feature of the ILLSS model is the remarkable stability of slender domains (tethers), both  in computations and experiments. This is in sharp contrast with computations that suggest singularity formation for tethers evolving by the Hele-Shaw flow~\cite{Almgren}. Our observation may be of value in addressing these harder questions.

\section{The model}
\label{sec:model}

The ILLSS model is as follows~\cite[\S 2.4]{AlBe}. The reader is referred to~\cite{AlBe} for a description of the assumptions that underlie the model.

The subfluid occupies the semi-infinite domain $B=\big\{(x,y,z)\in\R^3\big|z<0\big\}$.
The velocity $\vel = u \,\ei + v \,\ej + w \,\ek$  and pressure $P$ satisfy the Stokes equations
\ba
\label{eq:stokes2}
\nabla \cdot \vel =0, && z <0, \\
\label{eq:stokes1}
\lap \vel = \nabla P, && z <0.
\ea
The normal velocity $w$ vanishes on the boundary $z=0$. We assume that  $\vel$ and its derivatives decay sufficiently rapidly as $x^2+y^2+z^2 \to \infty$. 

The Langmuir layer is modeled as the surface $z=0$ of $B$. It decomposes into two 
{\it Langmuir layer domains} $\Omega$ and $\Omega^c$ separated by 
a moving boundary $\partial\Omega$. The motion of the boundary is determined by force balance and kinematic conditions. 

The force balance is as follows. The Langmuir layer $z=0$ is assumed to be in hydrostatic equilibrium.  Then the surface pressure $\Pi$ on $z=0$, the tangential surface stress $\Fs$, and the line tension $\Fl$ satisfy
\be\label{eq:balance}
\gradperp\Pi=\Fs+ \Fl,
\ee
where $\gradperp = \partial_x \, \ei + \partial_y \, \ej$.  The Stokes flow in the subfluid determines
\be\label{eq:tanstr}
\Fs(x,y)= -u_z(x,y,0) \, \ei - v_z (x,y,0) \, \ej. 
\ee
The surface pressure $\Pi$ jumps across the domain boundary $\partial \Omega$ because of the curvature and the line tension. This may be written
\be
\label{eq:linetension}
\Fl=\kappa\normal \delta(d)
\ee
where $\kappa$ is the curvature of $\partial \Omega$, $\normal$ is the unit outward normal, $\delta$ is the Dirac delta, and $d$ denotes the normal distance from $\partial \Omega$. The convention in~\cite{AlBe} is that $\kappa <0$ for convex $\Omega$.

We now consider the kinematic conditions. The surface velocity of the Langmuir layer is denoted $\Vel$. We assume the Langmuir layer is incompressible 
\be
\label{eq:vel_layer}
\gradperp \cdot \Vel =0,
\ee
and that the surface velocity is compatible with the subfluid velocity
\be
\label{eq:compatibility1}
\Vel(x,y) = \vel(x,y,0).
\ee
Finally, a point $\mathbf{\Gamma}$ on the boundary is advected by the surface velocity
\be
\label{eq:boundary}
\frac{D\mathbf{\Gamma}}{Dt} = \Vel.
\ee
This  completes the specification of the ILLSS model.

Alexander {\em et al} assume in addition that the flow field leaves horizontal sections invariant. As a consequence the vertical velocity $w=0$ vanishes in $B$ (not just on $z=0$). This assumption was made in earlier work~\cite{Lubensky,Stone} and allows considerable simplification via the introduction of streamfunctions.

\section{Gradient flow structure}
\label{sec:gradient}
\subsection{The framework}
The ingredients of a gradient flow structure are a manifold $\man$ with a metric $\metric$ and an energy functional $\energy$. 
If $m \in \man$, the gradient flow may be written in weak form as
\be
\label{eq:grad_weak} 
g(\dot{m},\vect) = -\langle  d\energy(m), \vect\rangle, \quad \vect \in T_m\man.
\ee
Here $\langle \cdot, \cdot \rangle$ denote the duality pairing between $1$-forms (elements of  $T_m\man^*$) and vectors (elements of $T_m \man$). The role of the metric is to convert the $1$-form $-d\energy(m)$ into a vector $\dot{m} \in T_m \man$. 

We now show that the ILLSS model is a gradient flow. The calculations are similar in spirit to~\cite{Otto}.
The manifold $\man$ consists of $C^\infty$ diffeomorphisms of $B$ that preserve orientation, volume and the surface area on the boundary $\partial B$. We denote a typical element of $\man$ by $\vpa$ and its restriction to $\partial B$ by $\vpb$. 
Diffeomorphisms in $\man$ correspond to flows generated by divergence free vector fields in $B$ with a divergence free trace on $\partial B$. Assume $\vell$ is
a $C^\infty$ vector field $\bar{B} \to \R^3$ with trace $\Vell(x,y)=\vell(x,y,0)$ and $\Vell \cdot \ek=0$. The associated flow $\vpa_\vell(\tau)$ is the solution to 
\be
\label{eq:flow_vel}
\partial_\tau \vpa_\vell(\tau) = \vell \circ \vpa_\vell(\tau), \quad \vpa(0)=\id.
\ee
This also implies
\be
\label{eq:flow_velb}
\partial_\tau \vpb_\vell(\tau) = \Vell \circ \vpb_\vell (\tau), \quad \vpb(0)=\id.
\ee
These calculations allow us to identify the tangent space $T_\vpa \man$. We first define the linear space of Eulerian velocity fields that form $T_\id \man$.
\be
\label{eq:euler_fields}
\tanv= \left\{ \vell \in C^\infty(\bar{B},\R^3) \left| \;\; \nabla \cdot \vell =0, \;\; \gradperp \cdot \Vell=0, \;\;\Vell\cdot \ek=0 \right. \right\}. 
\ee
Let $\vell_\vpa$ denote the composition $\vell \circ \vpa$. Then the tangent space at $\vpa \in \man$ is
\be
\label{eq:tangent_space}
T_\vpa \man = \left\{\vell_\vpa := \vell \circ \vpa,  \vell \in \tanv  \right\}. 
\ee 

We now introduce the metric. Let $\vpa \in \man$ and $\vell \in \tanv$. We define
\be
\label{eq:metric}
\metric_\vpa(\vell_\vpa,\vell_\vpa) = \int_B \tr \left( \nabla \vell^T \nabla \vell \right) \, dx = \sum_{i,j=1}^3 \int_B \left( \partial_{x_i} v_{j} \right)^2  \,dx.
\ee
If the integral vanishes, $\vell$ is a constant that vanishes because of the boundary condition. Thus, $\metric$ is indeed a metric. Note also that it is enough to specify the metric by  \qref{eq:metric}. This is equivalent to 
\be
\label{eq:metric2}
\metric_\vpa(\vell_\vpa,\velll_\vpa) = \int_B \tr \left( \nabla \vell^T \nabla \velll \right) \, dx, \quad \vell,\velll \in \tanv.
\ee
The principal heuristic idea is that the metric is determined by dissipation. All the dissipation in the ILLSS model is in the subfluid. The expression \qref{eq:metric} is precisely the dissipation in a Stokes flow in the subfluid. 

\subsection{The unconstrained ILLSS model as a gradient flow}
The natural energy of the ILLSS model is the perimeter of the interface $\partial \Omega_t$. If an initial phase configuration $\Omega_0 \subset \partial B$ is fixed, every $\vpa \in \man$ defines a new configuration $\vpb(\Omega_0)$ (recall that $\vpb$ is the restriction of $\vpa$ to $\partial B$). We set
\be
\label{eq:energy}
\energy(\vpa) = \Per(\vpb(\Omega_0)).
\ee

We are now ready to formulate the ILLSS model as a gradient flow. Equation \qref{eq:grad_weak} takes the form
\be
\label{eq:grad_Lang}
\metric_{\vpa}\left(\partial_t \vpa, \vell_\vpa \right) = -\langle d\energy(\vpa), \vell_\vpa \rangle, \quad \vell_\vpa \in T_\vpa \man, 
\ee
where the flow of diffeomorphisms $\vpa(t)$ is generated by the `true' velocity
\be
\label{eq:true_flow}
\partial_t \vpa = \vel \circ \vpa, \quad \partial_t \vpb = \Vel \circ \vpb,\quad \vpa(0)=\id, \;\;\vpb(0)=\id.
\ee
We show that equations \qref{eq:true_flow} and \qref{eq:grad_Lang} are equivalent to the free boundary evolution of  Section~\ref{sec:model}. 

We first compute the left hand side of \qref{eq:grad_Lang}. We use equation \qref{eq:metric2}, integrate by parts and use \qref{eq:tanstr} to obtain
\ba
\nonumber
\lefteqn{\metric_{\vpa}\left(\partial_t \vpa, \vell_\vpa\right) = \int_{B}\tr \left( \nabla\vel^T  \nabla \vell \right)\,dx}
\\
\nonumber
&& =-\int_{B}\lap \vel(x)\cdot \vell (x)\,dx
+\int_{\partial B}\vel_z\cdot \Vell \\
\label{eq:lhs_grad}
&&  =
-\int_{B}\lap \vel(x)\cdot\vell (x)\,dx -\int_{\partial B} \Fs \cdot\vell.
\ea

The right hand side of \qref{eq:grad_Lang} is computed as follows. The first variation of the perimeter is given by the curvature~\cite[\S 2.9]{Simon}. That is, 
\ba
\nonumber
\lefteqn{ \langle d\energy(\vpa), \vell \rangle = \left. \frac{d \energy \left(\vpa_\vell(\tau)\right)}{d\tau} \right|_{\tau=0}}\\
&& 
\label{eq:rhs_grad}
= -\int_{\partial \Omega}\kappa\normal \cdot\Vell=
-\int_{\partial B}\kappa\delta(d)\normal \cdot\Vell 
=-\int_{\partial B} \Fl \cdot\Vell.
\ea
(The $-$sign follows from the convention of~\cite{AlBe} that the curvature of a convex domain is negative). We combine \qref{eq:lhs_grad} and \qref{eq:rhs_grad} to see that \qref{eq:grad_Lang} is equivalent to
\be
\label{eq:weak_euler}
-\int_{B}\lap \vel(x)\cdot\vell (x)\,dx = \int_{\partial B} \left( \Fl +\Fs\right)  \cdot\Vell, \quad \vell \in \tanv,
\ee
with $\Fl$ and $\Fs$ defined by the constitutive relations \qref{eq:tanstr} and \qref{eq:linetension}. It follows that 
\be
\label{eq:weak1}
-\int_{B}\lap \vel(x)\cdot\vell (x)\,dx =0, 
\ee
for every $\vell \in \tanv$ with $\Vell=\mathbf{0}$. 

The Helmholtz decomposition for a vector field $\velf \in L^2(B)$ is
\be
\label{eq:helm1}
\velf = \nabla P + \velll, \quad \nabla \cdot \velll =0, \quad \Velll\cdot \ek=0.
\ee
The pressure $P$ is obtained as a solution to the Neumann problem
\be
\label{eq:helm2}
\lap P =0, \quad  z<0, \qquad \partial_z P = \velf \cdot \ek, \quad z=0.
\ee
$P$ is arbitrary up to a constant, but $\nabla P \in L^2(B)$ is unique, as is $\velll$. We apply this decomposition to $\velf=\lap u$. Then  \qref{eq:weak1} yields
\be
\label{eq:weak10}
-\int_{B}\velll(x)\cdot\vell (x)\,dx =0, 
\ee
for every $\vell \in \tanv$ with $\Vell=\mathbf{0}$. It follows that $\velll=\mathbf{0}$. Thus, $\vel$ solves Stokes equations \qref{eq:stokes2}--\qref{eq:stokes1}. Further, since $\lap \vel=\nabla P$, we integrate by parts in \qref{eq:weak_euler}, and use $\Vell\cdot k=0$ to find
\be
\label{eq:weak2}
\int_{\partial B} \left( \Fl +\Fs\right)  \cdot\Vell =0, 
\ee
for every smooth, divergence free vector field $\Vell$ on $\partial B$. We apply the Helmholtz decomposition again (but now in $\R^2$) to deduce the existence of a surface pressure $\Pi$ such that \qref{eq:balance} holds. Finally, the kinematic condition \qref{eq:boundary} is immediately implied by the fact that $\vpb(t)$
 solves \qref{eq:flow_velb} with $\Vell=\Vel$. This establishes the equivalence between \qref{eq:grad_Lang}--\qref{eq:true_flow} and the ILLSS model of Section~\ref{sec:model}.

\subsection{The constrained ILLSS model as a gradient flow}
We now consider the constrained flow considered in detail in~\cite{AlBe}. The constraint is that the flow leaves horizontal planes $z=$constant invariant. The flow is now restricted to the submanifold $\manb \subset \man$ that consists of $\vpa \in \man$ such that 
$\vpa(x,y,z) = \left( \vpa_1(x,y,z), \vpa_2(x,y,z),z\right)$. The corresponding tangent space $T_\vpa \manb$ consists of vector fields $\vell \circ \vpa$ as in \qref{eq:tangent_space} with the additional constraint $\vell \cdot \ek=0$. We now claim that the constrained ILLSS flow of~\cite{AlBe} may be written as the gradient flow
\be
\label{eq:grad_Langb}
\metric_{\vpa}\left(\partial_t \vpa, \vell_\vpa\right) = -\langle d\energy(\vpa), \vell_\vpa \rangle, \quad \vell_\vpa \in T_\vpa \manb, 
\ee
where $\vpa(t)$ also solves \qref{eq:true_flow}. Again, we must show that $\vel$ satisfies the equations of Section~\ref{sec:model}. 

The left and right hand sides of \qref{eq:grad_Lang} are computed exactly as earlier and yield \qref{eq:weak_euler} as earlier along with the constitutive relations \qref{eq:tanstr} and \qref{eq:linetension}. We again choose $\vell$ with $\Vell=0$ to  see that  \qref{eq:weak1} holds for every $\vell$ that is divergence free and satifies $\vell\cdot \ek =0$ and $\Vell=0$. The argument involving the Helmholtz decomposition of $\lap \vel$ is modified as follows. We now decompose a vector field $\velf \in L^2(B)$ with $\velf \cdot \ek=0$ in $B$ as in \qref{eq:helm1} with the additional constraint that $\velll \cdot \ek=0$ in $B$. It now follows that $P \equiv 0$, since the Neumann boundary condition in \qref{eq:helm2} is simply $\vel \cdot \ek=0$. We apply this decomposition with $\velf=\lap \vel$ to \qref{eq:weak1} and obtain \qref{eq:weak2} and the conclusion $\velll =\mathbf{0}$ as above. This shows that the velocity field $\vel$ is harmonic (thus also a solution of the Stokes equation with $P \equiv 0$).  It now follows immediately from \qref{eq:weak_euler} that \qref{eq:weak2} holds for every divergence free vector field $\Vell$. This yields \qref{eq:balance}. This shows the equivalence between \qref{eq:grad_Langb} and \qref{eq:true_flow} and the constrained model of \cite{AlBe}. The conclusion that $\vel$ is harmonic is deduced in~\cite{AlBe} using streamfunctions.

\section{Acknowledgements}
This work was supported by the National Science Foundation under grants DMS 05-30862 (MH) and  DMS 07-48482 (GM). 
We thank Andy Bernoff for encouraging our interest in this model.

\bibliographystyle{siam}
\bibliography{langmuir}

\end{document}